%% file: ModelOfMultiverse.tex
\documentclass{amsart}
\usepackage{diagrams,amssymb}
\diagramstyle[tight,centredisplay,dpi=600,noPostScript,heads=LaTeX]%
\input MathMacrosJDH
\newcommand{\x}{\mathfrak{X}}
\begin{document}
\author{Victoria Gitman}
\address{V. Gitman, New York City College of Technology of The City University of New York, Mathematics, 300 Jay Street, Brooklyn, NY 11201}
\email{vgitman@nylogic.org,
http://websupport1.citytech.cuny.edu/faculty/vgitman}
\author{Joel David Hamkins}
\address{J. D. Hamkins, Mathematics, The Graduate Center of The City University of New York, 365 Fifth Avenue, New York, NY 10016
 \& Mathematics, College of Staten Island of CUNY, Staten Island, NY 10314}
\email{jhamkins@gc.cuny.edu, http://jdh.hamkins.org}
\thanks{The research of the first author has been supported by grants from the
CUNY Research Foundation. The research of the second author has been
supported in part by grants from the CUNY
 Research Foundation and the National Science Foundation.}

\begin{abstract}
If\/ \ZFC\ is consistent, then the collection of countable
computably saturated models of \ZFC\ satisfies all of the Multiverse
Axioms of \cite{Hamkins:TheSet-TheoreticalMultiverse}.
\end{abstract}

\title{A natural model of the multiverse axioms}
\maketitle

\section{Introduction}

The multiverse axioms  that are the focus of this article arose in
connection with a continuing debate in the philosophy of set theory
between the \emph{Universe view}, which holds that there is a unique
absolute set-theoretical universe, serving as set-theoretic
background for all mathematical activity, and the \emph{Multiverse
view}, which holds that there are many set-theoretical worlds, each
instantiating its own concept of set. We refer the reader to
\cite{Hamkins:TheSet-TheoreticalMultiverse} and to several other
articles in the same special issue of the Review of Symbolic Logic
for a fuller discussion of this philosophical exchange (see also
\cite{Hamkins:TheMultiverse:ANaturalContext,
Hamkins2009:SomeSecondOrderSetTheory}). The multiverse axioms
express a certain degree of richness for the set-theoretic
multiverse, flowing from a perspective that denies an absolute
set-theoretic background.

Meanwhile, the multiverse axioms admit a purely
mathematical, non-philosophical treatment, on which we
shall focus here. We shall internalize the study of
multiverses to set theory by treating them as mathematical
objects within \ZFC, allowing for a mathematized simulacrum
inside $V$ of the full philosophical multiverse (which
would otherwise include universes outside $V$).
Specifically, in this article we define that a {\df
multiverse} is simply a nonempty set or class of models of
\ZFC\ set theory. The multiverse axioms then correspond to
the features listed in Definition
\ref{Definition.MultiverseAxioms}, which such a collection
may or may not exhibit.

\begin{definition}[Multiverse Axioms]
\rm Suppose that $\cal M$ is a multiverse, a nonempty
collection of models of
\ZFC.\label{Definition.MultiverseAxioms}
\begin{enumerate}
 \item The {\df Realizability} axiom holds for $\cal M$
     if whenever $M$ is a universe in $\cal M$ and $N$
     is a definable class of $M$ satisfying \ZFC\ from
     the perspective of $M$, then $N$ is in $\cal M$.
 \item The {\df Forcing Extension} axiom holds for
     $\cal M$ if whenever $M$ is a universe in $\cal M$
     and $\P$ is a forcing notion in $M$, then $\cal M$
     has a forcing extension of $M$ by $\P$, a model of
     the form $M[G]$, where $G$ is an $M$-generic
     filter for $\P$.
 \item The {\df Class Forcing Extension} axiom holds
     for $\cal M$ if whenever $M$ is a universe in
     $\cal M$ and $\P$ is a \ZFC-preserving class
     forcing notion in $M$, then $\cal M$ has a forcing
     extension of $M$ by $\P$, a model of the form
     $M[G]$, where $G$ is an $M$-generic filter for
     $\P$.
 \item The {\df Countability} axiom holds for $\cal M$
     if for every universe $M$ in $\cal M$ there is
     another universe $N$ in $\cal M$, such that $M$ is
     a countable set in $N$.
 \item The {\df Wellfoundedness Mirage} axiom holds for
     $\cal M$ if for every universe $M$ in $\cal M$,
     there is $N$ in  $\cal M$, which thinks $M$ is a
     set with an ill-founded $\omega$.
\end{enumerate}
Although the next axioms do not appear in
\cite{Hamkins:TheSet-TheoreticalMultiverse}, we shall
nevertheless consider them here. They follow a suggestion
of Jonas Reitz, who proposed that whenever a universe $M$
in the multiverse has a measurable cardinal, then it should
be the internal ultrapower of another universe $V$, sending
its critical point to that cardinal. That is, the
suggestion is that we should be able to iterate large
cardinal embeddings backwards. Here, we generalize the idea
to other ultrapowers and to embeddings generally.
\begin{enumerate}
 \item[(6)] The {\df Reverse Ultrapower Axiom} holds
     for $\cal M$ if for every universe $M$ in $\cal
     M$, there is a universe $N$ in $\cal M$ such that
     $M$ is the internal ultrapower of $N$ by an
     ultrafilter on $\omega$ in $N$.
 \item[(7)] The {\df Strong Reverse Ultrapower Axiom}
     holds for $\cal M$ if for every universe $M_1$ in
     $\cal M$ and every ultrafilter $U_1$ in $M_1$ on a
     set $X_1$ in $M_1$, there is $M_0$ in $\cal M$,
     with an ultrafilter $U_0$ on a set $X_0$, such
     that $M_1$ is the internal ultrapower of $M_0$ by
     $U_0$, sending $U_0$ to $U_1$.
\item[(8)] The {\df Reverse Embedding Axiom} holds for
    $\cal M$, if for every universe $M_1$ in $\cal M$
    and every embedding $j_1:M_1\to M_2$ definable in
    $M_1$ from parameters and thought by $M_1$ to be
    elementary, there is $M_0$ in $\cal M$ and similarly
    definable $j_0:M_0\to M_1$ in $M_0$, such that
    $j_1$ is the iterate of $j_0$, meaning
    $j_1=j_0(j_0)$.
\end{enumerate}
\end{definition}
In other words, the Reverse Embedding axiom asserts that every
internal elementary embedding $j_1:M_1\to M_2$ arises as an iterate
of an earlier embedding. The idea is that if we are living in $M_1$
and see the embedding $j_1$, then for all we know, it has already
been iterated an enormous number of times. To be precise, by
$j_1=j_0(j_0)$, we mean that if $j_0$ is definable by
$\varphi(x,y,a)$ over $M_0$, then $j_1$ is definable by
$\varphi(x,y,j_0(a))$ over $M_1$.

There is, of course, a certain degree of redundancy in the
axioms; for example, the Class Forcing Extension axiom
implies the Forcing Extension axiom, the Reverse Embedding
axiom implies the Reverse Ultrapower axioms, and the
Reverse Utrapower and Countability axioms imply the
Wellfoundedness Mirage. There are also a few subtler
points. In several of the axioms, when it is stated that
one model M is an element of another model N, what is meant
is that there is an object in $N$ that $N$ thinks is a pair
$\<m,E>$ for which $E$ is a binary relation on $m$, and
externally, the set $\set{a\in N\st N\satisfies a\in m}$
with the relation $\set{(a,b) \st N\satisfies aEb}$ is
isomorphic to $M$ with its relation. Another subtle issue
is that in the Countability axiom, although $M$ must be a
countable set in $N$, there is no insistence that $N$
regard $M$ as a model of \ZFC; indeed, since $N$ itself may
be nonstandard, it may have a nonstandard version of \ZFC,
with nonstandard size axioms that $M$ does not satisfy in
$N$, even if $M$ satisfies \ZFC\ externally. Similarly, the
Wellfoundedness Mirage axiom requires that $M$ is seen to
be ill-founded by $N$, but again $N$ may not look upon $M$
as a model of \ZFC, since $N$ may itself have nonstandard
size axioms, and there is no reason to expect that $N$
believes $M$ to satisfy them. The \emph{Realizability}
axiom, on the other hand, only applies to the definable
models of the universe that satisfy its own version of
\ZFC. For example, if $M$ is a model of \ZFC\ with
nonstandard $\omega$, then for every natural number $n$ in
$M$, there will be $V_\alpha^M$ that $M$ believes to model
the $\Sigma_n$ theory of its \ZFC, and when $n$ is
nonstandard, this includes all of the standard \ZFC; so
although these $V_\alpha^M$ are models of full \ZFC, the
Realizability axiom does not apply to them.

Hamkins \cite{Hamkins:TheSet-TheoreticalMultiverse}
provided a model of the first five multiverse axioms
constructed via iterated ultrapowers, and he inquired at
that time whether the collection of all countable
nonstandard models of set theory might already be a model
of the multiverse axioms. Observation
\ref{Observation.SaturationRequired} and Theorem
\ref{Theorem.Nonmodels} below show that this is too much,
since the axioms impose a requirement of computable
saturation on the models. Nevertheless, in this article we
prove that this is the only such obstacle, for our main
theorem shows that if \ZFC\ is consistent, then the
collection of all countable computably saturated models of
\ZFC\ satisfies all of the multiverse axioms.

\begin{maintheorem}\label{MainTheorem}
If\/ \ZFC\ is consistent, then the collection of all
countable computably saturated models of \ZFC\ satisfies
all the multiverse axioms.
\end{maintheorem}
We note that the main theorem will imply that it is not
true that every multiverse satisfying the axioms must
consist of countable models. For example, if $N$ is a
nonstandard model of $\ZFC+\Con(\ZFC)$, then we may let
$\cal M$ be the models that $N$ thinks are countable
computably saturated models of \ZFC. By the main theorem,
$N$ thinks that this collection satisfies all the
multiverse axioms, and it follows that this will really be
the case also outside of $N$. But the models in $\cal M$
will be at least as large as $\omega^N$, which could be as
large in cardinality as we like.

\section{Computably Saturated Models of Set Theory}

In this section, we shall explain precisely why we must
restrict to the computably saturated models and review the
key properties of these models that are needed for the
proof of the main theorem. Computable saturation was
introduced in \cite{schlipf:recursive_saturation} and is
also commonly known as \emph{recursive} saturation. A model
$M$ of a computable language $\mathcal L$ is said to be
\emph{computably saturated}, if for every finite tuple
$\bar{a}$ in $M$, every finitely realizable computable type
$p(\bar{a}, \bar{x})$ is already realized in $M$. A type
$p(\bar{y},\bar{x})$ in a computable language $\mathcal L$
is computable when the set of the G\"{o}del codes of its
formulas is a computable set in the usual sense of Turing
computability. In the future, we shall freely associate
types with subsets of $\N$ consisting of the G\"{o}del
codes of their formulas. A model of \ZFC\ set theory is
\emph{$\omega$-nonstandard} if it has a nonstandard
$\omega$. Because tuples can be viewed as a single set in
models of ZFC, for these it suffices to consider only
computable types of the form $p(a,x)$. Note that a
computably saturated model of \ZFC\ must necessarily be
$\omega$-nonstandard since the type $p(x)=\left\{{\bf
n}<x\mid n\in\omega\right\}\union\{x<\omega\}$, where ${\bf
n}$ is the term $1+\cdots+1$ with $n$ many $1$s, is a
finitely realizable computable type over any model of \ZFC.
For any model $M$ of \ZFC, we may consider the trace on the
natural numbers of the sets that exist in $M$.
Specifically, we say that a set $A\subseteq \N$ is
\emph{coded} by $a$ in $M$, if $a$ is a set in $M$ whose
intersection with the standard natural numbers is exactly
$A$. When $a$ is a set of natural numbers in $M$, then $A$
is also known as the \emph{standard part} of $a$, and the
collection of all sets $A$ arising this way is accordingly
called the \emph{Standard System} of $M$, denoted $SSy(M)$.
Note that the standard system of any model must include all
computable sets, since the model will agree on the behavior
of any computation that halts in the standard $\N$.
Standard systems have been extensively studied in the
context of models of Peano Arithmetic where they play a
crucial conceptual role, but the notion was originally
introduced for models of various set theories
\cite{friedman:ssy}.

\begin{observation}\label{Observation.SaturationRequired}
Any multiverse satisfying the Wellfoundedness Mirage axiom
must consist entirely of computably saturated models of
\ZFC.
\end{observation}

\begin{proof}
If $\cal M$ satisfies the Wellfoundedness Mirage axiom, it
follows that every member of $\cal M$ has a nonstandard
$\omega$. Thus, again by the Wellfoundedness Mirage axiom,
every member of $\cal M$ is a set in a model of \ZFC\
having a nonstandard $\omega$. Thus, every member of $\cal
M$ is computably saturated by Lemma \ref{Lemma.Saturation}.
\end{proof}

\begin{lemma}\label{Lemma.Saturation}
Every model of \ZFC\ that is an element of an
$\omega$-nonstandard model of \ZFC\ is computably
saturated.
\end{lemma}

\begin{proof}
Suppose that $M$ is a model of \ZFC\ that is an element of
an $\omega$-nonstandard model $N$ of \ZFC. In order to see
that $M$ is computably saturated, suppose that $p(b,x)$ is
a computable finitely realizable type over $M$. Let $a\in
N$ code $p(y,x)$. Since $p(b,x)$ is finitely realizable and
$N$ has a truth predicate for $M$, for every $n\in \N$, the
model $N$ knows that there is $c\in M$ such that
$M\models\varphi(b,c)$ for every formula $\varphi(y,x)$
with G\"{o}del code less than $n$ in $a$. Because the
standard $\N$ is not definable in $N$, there must be a
nonstandard natural number $d\in N$ and an element $e\in M$
with $N$ satisfying that $M\models \varphi(b,e)$ for every
(possibly nonstandard) formula $\varphi(y,x)$ with
G\"{o}del code less than $d$ in $a$. Since $d$ is
nonstandard, this includes every formula in $p(y,x)$. By
the absoluteness of satisfaction for standard formulas, it
follows that $M\models\varphi(b,e)$ for every
$\varphi(b,x)$ in $p(b,x)$ and thus $e$ realizes $p(b,x)$.
\end{proof}

\begin{lemma}\label{Lemma.Existence}
If\/ \ZFC\ is consistent, then there are $2^{\aleph_0}$
many pairwise non-isomorphic countable computably saturated
models of \ZFC. Every real is in the standard system of
such a model.
\end{lemma}

\begin{proof}
If \ZFC\ consistent, then every completion of \ZFC\ as a
theory has a countable computably saturated model, because
any countable model of \ZFC\ can be extended elementarily
to a computably saturated model by successively realizing
types in a countable elementary chain. For any real $x$,
one can ensure that the type expressing that $x$ is coded
is realized.
\end{proof}
A model $M$ of \ZFC\ is said to be $SSy(M)$-saturated if it
realizes every finitely realizable type coded in $M$. It
turns out that a model $M$ is computably saturated if and
only if it is $SSy(M)$-saturated. To see this, fix a type
$p(y,x)$ coded by $a\in M$ and fix $b\in M$ such that
$p(b,x)$ is finitely realizable. Define a new type
$q(b,a,x)$ to consist of all formulas of the form
$(\ulcorner\varphi(y,x)\urcorner\in a)\rightarrow
\varphi(b,x)$, and observe that $q(y,z,x)$ is computable
and finitely realizable, using the objects realizing the
corresponding fragment of $p(b,x)$. Thus, there is some $e$
in $M$ realizing $q(b,a,e)$, and it follows that $e$
realizes $p(b,e)$, as desired. We note also that the type
of any element in a computably saturated model is in the
standard system of that model: for $a\in M$, define
$p(a,x)$ to be the type consisting of all formulas of the
form $\ulcorner \varphi(y)\urcorner\in x\leftrightarrow
\varphi(a)$ and observe that it is computable and finitely
realizable; thus, the type of $a$ is coded in $M$. In
particular, the theory of any computably saturated model is
an element of its standard system.  According to
\cite{kaye:modelsofpa}, $SSy(M)$-saturation was introduced
by Wilmers in his unpublished 1975 thesis where he
established the above equivalence.

The next lemma generalizes another fundamental result from models of
PA that appears in \cite{smorynsky:recsat.characterization} but has
as well been attributed to Jensen and Ehrenfeucht
\cite{ehrenfeucht.recsat}, and Wilmers, among others.
\begin{keylemma}\label{Lemma.Isomorphism}
Any two countable computably saturated models of \ZFC\ with
the same theory and the same standard system are
isomorphic.
\end{keylemma}

\begin{proof}
This is a standard model-theoretic back-and-forth construction. The
observations above ensure that the models are standard
system-saturated, and all types of their elements are coded in the
standard system. Thus, we may construct the desired isomorphism in a
countable recursive procedure that maps elements of one model to
elements in the other realizing the same types over what has been
defined so far.
\end{proof}

The following lemma will be critical for our verification
of the Wellfoundedness Mirage axiom in the Main Theorem.
This fact may have been known some time ago. For example,
Schlipf
\cite[III.2.6]{Schlipf1977:AGuideToIdentificationOfAdmissibleSets:}
proved that every computably saturated model of \ZF\ is an
element of an $\omega$-nonstandard model of \ZF, and
Ressayre \cite[3.3]{Ressayre1983:ModelesNonStandard} proved
that every model of \ZF\ is elementarily equivalent to a
model of \ZF\ containing as an element an isomorphic copy
of itself. (See \cite{Halimi:ModelsAndUniverses} for an
interesting discussion.)

\begin{lemma}\label{Lemma.ContainsACopy}
Every countable computably saturated model of \ZFC\
contains an isomorphic copy of itself as an element, which
it thinks is $\omega$-nonstandard. That is, if $M$ is a
countable computably saturated model of \ZFC, then $M$ has
an element $N$ which it thinks is a countable
$\omega$-nonstandard model of a fragment of set theory,
such that $M\iso N$.
\end{lemma}

\begin{proof} Suppose that M is a countable
computably saturated model of ZFC. As we noted above,
$\Th(M)$ is coded by some $a\in M$. By the Reflection
Theorem, every finite subset of this theory is true in some
rank initial segment of $M$, and $M$ recognizes this for
any particular such finite subset. Since the standard cut
$\N$ is not definable in $M$, there must be a nonstandard
natural number $b$ in $M$, such that $M$ thinks the theory
consisting of all formulas whose G\"{o}del codes are in $a$
and less than $b$ is consistent. Since $b$ is nonstandard,
this includes the entire $\Th(M)$. By the Completeness
Theorem in $M$, therefore, we may build a model $N$ in $M$
satisfying this consistent fragment of $a$, which includes
all of $\Th(M)$, such that additionally, $M$ thinks $N$ is
$\omega$-nonstandard. Since $\omega^M$ is an initial
segment of $\omega^N$ and $M$ is $\omega$-nonstandard, it
follows that $M$ and $N$ have the same standard system.
Also, since $M$ is $\omega$-nonstandard, it follows by
Lemma \ref{Lemma.Saturation} that $N$ is computably
saturated. We conclude by Lemma \ref{Lemma.Isomorphism}
that actually $M\iso N$.
\end{proof}
By considering the situation from the perspective of the
smaller copy of the model, we deduce:
\begin{corollary}\label{cor:element}
Every countable computably saturated model of \ZFC\ is an
element of another countable computably saturated model of
\ZFC\ that thinks it is a countable $\omega$-nonstandard
model of a (nonstandard) fragment of set theory.
\end{corollary}
Note in Lemma \ref{Lemma.ContainsACopy} that although we know on the
outside that $N\satisfies\ZFC$, since it satisfies $\Th(M)$, it
could happen that $M\not\satisfies$``$N\satisfies\ZFC$,'' since
perhaps $M$ thinks that some of the nonstandard \ZFC\ axioms of $M$
fail in $N$. Despite this, Corollary \ref{cor:element} suffices to
verify the Countability and Well-foundedness Mirage axioms for the
collection of countable computably saturated models of \ZFC, since
as we mentioned there was no insistence in the axioms that the
larger model look upon the smaller as a model of what it thinks is
full \ZFC. Nevertheless, under a stronger assumption it is possible
to obtain the stronger conclusion. Surely a stronger assumption is
required, since if $N\satisfies\ZFC+$``$M\satisfies\ZFC$'', then
$N\satisfies\ZFC+\Con(\ZFC)$, and so $\Con(\ZFC+\Con(\ZFC))$. And if
this $N$ is an element of a further such model, then we get
$\Con(\Con(\Con(\ZFC)))$, and so on transfinitely. The stronger
assumption we shall make is that for every countable computably
saturated model $M$ of ZFC, the theory
$T_M=\ZFC+\{\Con(\ZFC+\Gamma)\mid \Gamma\subseteq_{\text{Fin}}
\Th(M)\}$ is consistent.

\begin{theorem} If\/ $M$ is a computably saturated countable model
of\/ \ZFC, then there is a countable computably saturated
model $N$ of\/ \ZFC\ containing $M$ as an element and
satisfying that $M$ is a nonstandard model of\/ \ZFC\ if
and only if the theory $T_M$ is consistent.
\end{theorem}
\begin{proof}
The forward implication is immediate, since any such model
$N$ will satisfy the theory $T_M$. For the converse
implication, suppose that $T_M$ is consistent. By Lemma
\ref{Lemma.Isomorphism}, it suffices to show that there
exists $N$ containing a countable computably saturated
model $K$ with the same theory and standard system as $M$
that it recognizes as a model of \ZFC. It will immediately
follow that $SSy(N)=SSy(K)$ and so we shall need to ensure
that $SSy(N)=SSy(M)$. Scott observed in
\cite{Scott1962:AlgebrasOfSets} that every standard system
is a \emph{Scott set}, that is, a Boolean algebra of
subsets of natural numbers that is closed under relative
computability and contain at least one branch through every
element that is a binary tree. In that paper, he famously
showed that given a countable Scott set $\x$ and a theory
$T\in \x$ extending PA, there is a model of $T$ whose
standard system is exactly $\x$. Wilmers in his thesis,
showed that this easily generalizes to obtaining a
computably saturated model. Let $\x=SSy(M)$ and observe
that $T_M\in \x$, since $T_M$ is computable in $\Th(M)$,
which is an element of $\x$. Summarizing, we can build a
countable computably saturated model $N$ of \ZFC\
satisfying $T_M$ and having same standard system as $M$.
Since $N$ satisfies $T_M$, it satisfies $\Con(\ZFC+\Gamma$)
where $\Gamma$ is a nonstandard segment containing
$\Th(M)$. So $N$ can build a countable model $K$ that it
thinks satisfies $\ZFC+\Gamma$.
\end{proof}

We have observed that the assumption that $T_M$ is
consistent transcends $\Con(\ZFC)$. But the assumption is
not so strong, for if $M$ is an element of an
$\omega$-model $N$ of \ZFC, then $N$ satisfies $T_M$. In
particular, if there is a transitive model $N$ of \ZFC,
then it satisfies $T_M$, and hence also $\Con(T_M)$, for
every countable model $M$ in $N$.

Let us close this section by mentioning the concept of
\emph{resplendency}, a powerful generalization of computable
saturation that has unified many applications of it. Resplendency is
a second-order analogue of computable saturation, in that it
concerns realizing second-order types; that is, it is about
interpreting a new predicate symbol on the universe. Specifically, a
first order structure $M$ is \emph{resplendent}, if every
finitely-realized computable type $p(X,\vec a)$ in the language of
$M$ expanded by a predicate symbol $X$ with $\vec a$ a finite list
of parameters from $M$ is realized in $\<M,X>$ for some
interpretation of $X$. (The type is finitely realized if all finite
subsets of $p$ are realized in such a model $\<M,X>$.) The concept
of resplendency was introduced by Barwise and Schlipfe
\cite{BarwiseSchlipf1976:IntroRecursivelySaturatedAndResplendentModels},
and independently by Ressayre
\cite{Ressayre1983:IntroductionToRecursivelySaturatedModels}, who
proved that every countable computably saturated model is
resplendent (see also
\cite{Smorynski1981:RecursivelySaturatedNonstandardModelsOfArithmetic}).
Schlipfe \cite{Schlipf1980:RecursivelySaturatedModelsOfSetTheory}
proved that a countable model of set theory is computably saturated
if and only if it is $\omega$-nonstandard and there is a club of
ordinals $\alpha$ with $V_\alpha\elesub V$. Moschovakis and Chang
(see \cite{ChangKeisler1990:ModelTheory}) proved that every
saturated model is resplendent. Although we have presented our
arguments in an elementary manner appealing only to computable
saturation, it appears that many of our lemmas can be fruitfully
generalized, by proving them via resplendency.

\section{Proof of Main Theorem}

Let us now complete the proof of the main theorem, which we
restate here for convenience.

\newtheorem*{maintheorem*}{Main Theorem}
\begin{maintheorem*}
If\/ \ZFC\ is consistent, then the collection of all
countable computably saturated models of \ZFC\ satisfies
all the multiverse axioms.
\end{maintheorem*}

\begin{proof}
We shall argue in turn that the collection $\cal M$ of all
countable computably saturated models of \ZFC\ satisfies
each of the multiverse axioms. First, since we have assumed
that \ZFC\ is consistent, Lemma \ref{Lemma.Existence} shows
that in fact there are many countable computably saturated
models of \ZFC. So $\cal M$ is nonempty.

Consider now the Realizability axiom. Suppose that
$M\in\cal M$ and $N$ is a definable class in $M$ and a
model of \ZFC. Since $M$ is an element of some other
nonstandard model $M'$ by Corollary \ref{cor:element}, it
follows that $N$ is also an element of $M'$, and so by
Lemma \ref{Lemma.Saturation}, it follows that $N$ is
computably saturated. Since $N$ is clearly also countable,
as $M$ was countable, it follows that $N\in\cal M$. Thus,
$\cal M$ satisfies the Realizability axiom.

For the Forcing axioms, suppose that $M\in\cal M$ and $\P$
is a forcing notion in $M$. Certainly we can easily produce
by diagonalization an $M$-generic filter $G\of\P$ and form
the forcing extension $M[G]$. Furthermore by Corollary
\ref{cor:element}, we can do so inside any model $M'$ which
looks upon $M'$ as countable. Thus, there is a forcing
extension $M[G]$ inside such an $M'$. It now follows by
Lemma \ref{Lemma.Saturation} that $M[G]$ is computably
saturated, as desired. So $\cal M$ satisfies the Forcing
and Class Forcing Extension axioms.

The difficult cases of the Wellfoundedness Mirage and
Countability axioms are exactly provided for by Corollary
\ref{cor:element}.

The Reverse Ultrapower axioms follow from the Reverse Embedding
axiom, so it suffices to consider that axiom. Suppose that $M_1$ is
countable and computably saturated and $j_1:M_1\to M_2$ is an
elementary embedding in $M_1$, defined in $M_1$ from some parameter
$z$, so that $j_1(x)=y\iff M_1\satisfies\varphi(x,y,z)$. By
interpreting this definition in $M_2$ using $j_1(z)$ we obtain the
iterate embedding $j_2=j_1(j_1):M_2\to M_3$, defined by
$j_2(x)=y\iff M_2\satisfies\varphi(x,y,j_1(z))$. Since the critical
point of $j_1$ must be at least $\omega^{M_1}$, which is
nonstandard, it follows that $M_1$ and $M_2$ share a nonstandard
initial segment of their natural numbers and therefore have the same
standard system. Since they also have the same theory, it follows by
Lemma \ref{Lemma.Isomorphism} that there is an isomorphism
$\pi:M_1\iso M_2$. Since the type of $z$ in $M_1$ is the same as the
type of $j_1(z)$ in $M_2$, we may assume in the back-and-forth
argument that $\pi(z)=j_1(z)$. Thus, since $j_1$ is defined in $M_1$
by $\varphi(x,y,z)$, the map $\pi$ carries $j_1$ to the class
defined in $M_2$ by $\varphi(x,y,\pi(z))$, which is $j_2$. In other
words, $\pi$ carries the entire map $j_1:M_1\to M_2$ isomorphically
to the map $j_2:M_2\to M_3$. And since $j_2=j_1(j_1)$ is by
definition an iterate of $j_1$, the Reverse Embedding axiom holds in
the case of $j_2:M_2\to M_3$. Since this is isomorphic via $\pi$ to
$j_1:M_1\to M_2$, it follows by replacing the objects with their
image under $\pi$ that there is $j_0:M_0\to M_1$ such that
$j_1=j_0(j_0)$, as desired.
\end{proof}

Recall that a model $M$ is said to be $\kappa$-\emph{saturated} for
a cardinal $\kappa$ if every finitely realizable type in the
language extended to include some $<\kappa$-many constants for
elements of the model is already realized. A  model of cardinality
$\kappa$ is said to be simply \emph{saturated} if it is
$\kappa$-saturated. It is a basic fact that any two saturated models
of the same theory and same cardinality are isomorphic. If $M$ is a
saturated model of ZFC and $N$ is a model of ZFC that is an element
of $M$, then $N$ must be saturated and have the same cardinality as
$M$. The cardinality is the same since $\omega^M$, by saturation, is
already of the same cardinality as $M$. For details on saturated
models, see \cite{ChangKeisler1990:ModelTheory}. Thus, it easily
follows that every saturated model of ZFC of cardinality $\kappa$
has an isomorphic copy of itself that it thinks is a countable
$\omega$-nonstandard model of a finite fragment of ZFC. Other facts
necessary for the proof of the Main Theorem follow for saturated
models of ZFC of cardinality $\kappa$ as well; in most cases they
are easier to see than for computable saturation because any two
elementarily equivalent saturated models of the same cardinality are
isomorphic. Thus, we get the following corollary of the Main
Theorem.
\begin{corollary}
If there are saturated models of \ZFC\ of cardinality $\kappa$, then
the collection of these satisfies all the multiverse axioms.
\end{corollary}

It is natural to wonder whether the collection of all
models of \ZFC\ forms a model of the multiverse axioms, or
whether the collection of all countable models of \ZFC\
does so. Unfortunately, neither does.

\begin{theorem}\label{Theorem.Nonmodels}
If\/ \ZFC\ is consistent, then the collection of all models
of\/ \ZFC\ is not a model of the multiverse axioms. Neither
is the collection of all countable models of\/ \ZFC, nor
the collection of all countable nonstandard models of\/
\ZFC, nor the collection of countable $\omega$-nonstandard
models of \ZFC, nor the collection of such models
restricted to a given consistent completion of \ZFC.
\end{theorem}

\begin{proof}
By Observation \ref{Observation.SaturationRequired}, all we
need to do for the first part is to show that there is a
model of \ZFC\ that is not computably saturated. In fact,
every consistent completion of \ZFC\ has a countable
$\omega$-nonstandard model that is not computably saturated
(and this proves the subsequent claims). To see this, take
any countable nonstandard model $M$ of \ZFC. The {\df
definable cut} of $M$ consists of all $x\in M$ such that
$x\in (V_\alpha)^M$, where $\alpha$ is a definable ordinal
in $M$ (without parameters). If $M_0$ is the definable cut
of $M$, then it is relatively easy to verify the
Tarski-Vaught criterion, and so $M_0\elesub M$. It follows
that $M_0$ has exactly the same definable ordinals as $M$,
and these are unbounded in the ordinals of $M_0$. Thus,
$M_0$ omits the type $p(x)$ asserting that whenever there
is a unique ordinal satisfying $\varphi(y)$, then $y<x$.
This is a computable finitely realizable type not realized
in $M_0$, and so $M_0$ is not computably saturated. Thus,
by Observation \ref{Observation.SaturationRequired}, it can
not be in any model of the multiverse axioms.
\end{proof}

Let us conclude this paper by considering the degree to which we
might expect a multiverse to be upward directed. Specifically, a
multiverse $\cal M$ is \emph{upward directed}, if for any two
elements $M,N\in{\cal M}$ there is an element $W\in{\cal M}$
containing (isomorphic copies of) $M$ and $N$ as elements. The
multiverse $\cal M$ is \emph{countably upward directed} if for any
countable subcollection ${\cal M}_0\of{\cal M}$, there is an element
$W\in{\cal M}$ containing (an isomorphic copy of) every element of
${\cal M}_0$. It is easy to see that the multiverse of all countable
computably saturated models of \ZFC\ is not upward directed. This is
because any two elements of an upward directed multiverse $\cal M$
containing only $\omega$-nonstandard models must have the same
standard system. Suppose that $M$ and $N$ are elements of an upward
directed multiverse $\cal M$ containing only $\omega$-nonstandard
models. By directedness, there is $W\in{\cal M}$ with $M$ and $N$
both in $W$. Since the $\omega^W$ is an initial segment of
$\omega^M$ and $\omega^N$, and is itself nonstandard, it follows
that all three models $M$, $N$ and $W$ have the same standard
system. Thus, all models in $\cal M$ have the same standard system.
Since any real can be placed into the standard system of a countable
computably saturated model of \ZFC, it follows that not all
countable computably saturated models of \ZFC\ have the same
standard system. So this multiverse is not upward directed.
Nevertheless, this is the only obstacle.

\begin{theorem}
If\/ \ZFC\ is consistent, the multiverse of countable computably
saturated models having a fixed standard system is countably upward
directed, and continues to satisfy all the multiverse axioms.
\end{theorem}

\begin{proof}
Fix a given Scott set $S$ and consider the multiverse ${\cal M}_S$
of all countable computably saturated models of \ZFC\ having
standard system $S$. We observe first that the proof of the main
theorem goes through for ${\cal M}_S$, since in each part of that
argument, the desired universe had the same standard system as the
original model. So it remains only to argue that ${\cal M}_S$ is
countably upward directed. Suppose that ${\cal
M}_0=\set{M_0,M_1,\ldots}$ is a countable subcollection of ${\cal
M}_S$, so that every $M_n$ is a countable computably saturated model
of \ZFC\ with standard system $S$. By the remarks before Lemma
\ref{Lemma.Isomorphism}, it follows that $\Th(M_n)$ is in $S$ for
every $n$. Let $M$ be any $\omega$-nonstandard model having standard
system $S$. Since $\Th(M_n)$ is coded in $M$, by arguments of the
proof of Lemma \ref{Lemma.ContainsACopy}, $M$ can build a model
$m_n$ satisfying the theory $\Th(M_n)$ and having $SSy(m_n)=S$.
Therefore, by Lemma \ref{Lemma.Isomorphism}, it follows that $m_n$
and $M_n$ are isomorphic. In summary, we have proved that {\it
every} model in ${\cal M}_S$ serves as a witness to the countable
upward directedness of ${\cal M}_S$.
\end{proof}

\bibliographystyle{alpha}
\bibliography{HamkinsBiblio,MathBiblio,multiverse}

\end{document}

%% file: MathMacrosJDH.tex
\newtheorem{theorem}{Theorem}
\newtheorem{maintheorem}[theorem]{Main Theorem}
\newtheorem{corollary}[theorem]{Corollary}

\newtheorem{lemma}[theorem]{Lemma}
\newtheorem{keylemma}[theorem]{Key Lemma}

\newtheorem{observation}[theorem]{Observation}

\newtheorem{definition}[theorem]{Definition}

\newcommand{\QED}{\end{proof}}

\def\BF#1.{{\bf #1.}}
\newcommand{\url}[1]{{\tt #1}}
\newcommand{\cal}{\mathcal}

\newcommand{\N}{{\mathbb N}}
\renewcommand{\P}{{\mathbb P}}

\newcommand{\of}{\subseteq}

\newcommand{\set}[1]{\{\,{#1}\,\}}

\newcommand{\elesub}{\prec}

\newcommand{\Th}{\mathop{\rm Th}}
\newcommand{\Con}{\mathop{{\rm Con}}}

\newcommand{\satisfies}{\models}

\newcommand{\union}{\cup}

\newcommand{\smalllt}{\mathrel{\mathchoice{\raise2pt\hbox{$\scriptstyle<$}}{\raise1pt\hbox{$\scriptstyle<$}}{\raise0pt\hbox{$\scriptscriptstyle<$}}{\scriptscriptstyle<}}}
\newcommand{\smallleq}{\mathrel{\mathchoice{\raise2pt\hbox{$\scriptstyle\leq$}}{\raise1pt\hbox{$\scriptstyle\leq$}}{\raise1pt\hbox{$\scriptscriptstyle\leq$}}{\scriptscriptstyle\leq}}}

\newcommand{\boolval}[1]{\mathopen{\lbrack\!\lbrack}\,#1\,\mathclose{\rbrack\!\rbrack}}
\def\[#1]{\boolval{#1}}

\newcommand{\UnderTilde}[1]{{\setbox1=\hbox{$#1$}\baselineskip=0pt\vtop{\hbox{$#1$}\hbox to\wd1{\hfil$\sim$\hfil}}}{}}
\newcommand{\Undertilde}[1]{{\setbox1=\hbox{$#1$}\baselineskip=0pt\vtop{\hbox{$#1$}\hbox to\wd1{\hfil$\scriptstyle\sim$\hfil}}}{}}
\newcommand{\undertilde}[1]{{\setbox1=\hbox{$#1$}\baselineskip=0pt\vtop{\hbox{$#1$}\hbox to\wd1{\hfil$\scriptscriptstyle\sim$\hfil}}}{}}
\newcommand{\UnderdTilde}[1]{{\setbox1=\hbox{$#1$}\baselineskip=0pt\vtop{\hbox{$#1$}\hbox to\wd1{\hfil$\approx$\hfil}}}{}}
\newcommand{\Underdtilde}[1]{{\setbox1=\hbox{$#1$}\baselineskip=0pt\vtop{\hbox{$#1$}\hbox to\wd1{\hfil\scriptsize$\approx$\hfil}}}{}}

\newcommand{\st}{\mid}

\newcommand{\iso}{\cong}
\def\<#1>{\langle#1\rangle}

\newcommand{\ZFC}{{\rm ZFC}}
\newcommand{\ZF}{{\rm ZF}}

\newcommand{\cell}[1]{\boxit{\hbox to 17pt{\strut\hfil$#1$\hfil}}}
\newcommand{\head}[2]{\lower2pt\vbox{\hbox{\strut\footnotesize\it\hskip3pt#2}\boxit{\cell#1}}}
\newcommand{\boxit}[1]{\setbox4=\hbox{\kern2pt#1\kern2pt}\hbox{\vrule\vbox{\hrule\kern2pt\box4\kern2pt\hrule}\vrule}}
\newcommand{\Col}[3]{\hbox{\vbox{\baselineskip=0pt\parskip=0pt\cell#1\cell#2\cell#3}}}
\newcommand{\tapenames}{\raise 5pt\vbox to .7in{\hbox to .8in{\it\hfill input: \strut}\vfill\hbox to
.8in{\it\hfill scratch: \strut}\vfill\hbox to .8in{\it\hfill output: \strut}}}
\newcommand{\Head}[4]{\lower2pt\vbox{\hbox to25pt{\strut\footnotesize\it\hfill#4\hfill}\boxit{\Col#1#2#3}}}
\newcommand{\Dots}{\raise 5pt\vbox to .7in{\hbox{\ $\cdots$\strut}\vfill\hbox{\ $\cdots$\strut}\vfill\hbox{\
$\cdots$\strut}}}
\newcommand{\df}{\it} 